\documentclass[a4paper,english,10pt,oneside]{article}\def\zibreport{1}

\usepackage{ifthen}
\usepackage{graphicx}
\usepackage{multicol}
\usepackage{footmisc}

\usepackage{ifthen}
\usepackage{amsmath,amsfonts,amssymb}
\usepackage{tabularx,supertabular,booktabs}
\usepackage[textsize=scriptsize, textwidth=3.3cm]{todonotes}
\usepackage{xspace}
\usepackage{lscape}
\usepackage{multirow}
\usepackage{tabularx}
\usepackage{comment}
\usepackage[utf8]{inputenc}

\usepackage{tikz}
\usetikzlibrary{matrix,calc,positioning}
\usetikzlibrary{petri}	
\usetikzlibrary{shapes.misc}
\tikzset{cross/.style={cross out, draw=black, minimum size=2*(#1-\pgflinewidth), inner sep=0pt, outer sep=0pt},
	cross/.default={1pt}}
\usepackage{subcaption}
\usepackage{pgfplots}

\usepackage{url}

\usepackage[bookmarks]{hyperref}
\hypersetup{
  pdftitle=First Experiments with Structure-Aware Presolving for a Parallel Interior-Point Method,
  pdfauthor={Ambros Gleixner, Nils-Christian Kempke, Thorsten Koch, Daniel Rehfeldt, Svenja Uslu}
}

\newcommand{\myorcidlink}[1]{\,\href{https://orcid.org/#1}{\raisebox{-0.45ex}{\includegraphics[width=1.8ex]{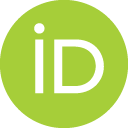}}}}
\newcommand{\myurl}[1]{\textsf{\footnotesize \url{#1}}\xspace}

\newcommand{\allcaps}[1]{\protect\scalebox{0.93}{#1}}

\newcommand{\LP}{\allcaps{LP}\xspace}
\newcommand{\LPs}{\allcaps{LPs}\xspace}

\newcommand{\MPI}{\allcaps{MPI}\xspace}
\newcommand{\HPC}{\allcaps{HPC}\xspace}

\newcommand{\R}{\mathbb{R}}

\newcommand{\solver}[1]{\textsc{#1}\xspace}

\newcommand{\soplex}{\solver{\allcaps{SoPlex}}}

\newcommand{\pips}{\solver{\allcaps{PIPS-IPM}}}

\ifthenelse{\zibreport = 1}{
	
	\usepackage{zibtitlepage}

	\ZTPAuthor{Ambros Gleixner\protect\myorcidlink{0000-0003-0391-5903}, Nils-Christian Kempke\protect\myorcidlink{0000-0003-4492-9818}, Thorsten Koch\protect\myorcidlink{0000-0002-1967-0077}, Daniel Rehfeldt\protect\myorcidlink{0000-0002-2877-074X} and Svenja Uslu}
	\ZTPTitle{\bf First Experiments with Structure-Aware Presolving for a Parallel Interior-Point Method}
	\ZTPInfo{}
	\ZTPNumber{19-39}
	\ZTPMonth{July}
	\ZTPYear{2019}
	
	\title{\LARGE\bf First Experiments with Structure-Aware Presolving for a Parallel Interior-Point Method}
	
	\author{
		Ambros Gleixner\protect\myorcidlink{0000-0003-0391-5903}, Nils-Christian Kempke\protect\myorcidlink{0000-0003-4492-9818}, Thorsten Koch\protect\myorcidlink{0000-0002-1967-0077},\\ Daniel Rehfeldt\protect\myorcidlink{0000-0002-2877-074X} and Svenja Uslu\\[0.5ex]
		\small
		Zuse Institute Berlin, Department of Mathematical Optimization,\\[-0.1ex]
		\small
		\texttt{\{gleixner,kempke,koch,rehfeldt\}@zib.de}\\[3ex]
	}
	
	\date{\normalsize July 26, 2019}	
	
}{
	\title{First Experiments with Structure-Aware Presolving for a Parallel Interior-Point Method}	
	\author{Ambros Gleixner\inst{1}\and Nils-Christian Kempke\inst{1}\and Thorsten Koch\inst{1,2}\and Daniel Rehfeldt\inst{1,2}\and Svenja Uslu\inst{1}}
	\institute{Zuse Institute Berlin, Takustr.~7, 14195~Berlin, Germany\\%
		\and Technische Universit{\"a}t Berlin, Stra{\ss}e des 17.~Juni~136, 10623~Berlin}
	\titlerunning{Structure-Aware Presolving for a Parallel Interior-Point Method}
	\authorrunning{A. Gleixner et al.}
}

\begin{document}
	
\ifthenelse{\zibreport = 1}{
	\zibtitlepage
}{\mainmatter}

\maketitle

\begin{abstract}	
In linear optimization, matrix structure can often be exploited algorithmically.
However, beneficial presolving reductions sometimes destroy the special structure of a given problem.
In this article, we discuss structure-aware implementations of presolving as part of a parallel interior-point method to solve linear programs with block-diagonal structure, including both linking variables and linking constraints.
While presolving reductions are often mathematically simple, their implementation in a high-performance computing environment is a complex endeavor.
We report results on impact, performance, and scalability of the resulting presolving routines on real-world energy system models with up to 700~million nonzero entries in the constraint matrix.
\ifthenelse{\zibreport = 0}{\keywords{block structure, energy system models, \HPC, linear programming, interior-point methods, parallelization, presolving}}{}

\end{abstract}	

\section{Introduction}
\label{sec:intro}

Linear programs (\LPs) from energy system modeling and from other applications based on time-indexed decision variables often exhibit a distinct block-diagonal structure.
Our extension~\cite{PIPSgor17} of the parallel interior-point solver \pips~\cite{pipsPetra} exploits this structure even when both linking variables and linking constraints are present simultaneously.
It was designed to run on high-performance computing (\HPC) platforms to make use of their massive parallel capabilities.
In this article, we present a set of highly parallel presolving techniques that improve \pips's performance while preserving the necessary structure of a given \LP.
We give insight into the implementation and the design of said routines and report results on their performance and scalability.

\begin{figure}\small\begin{alignat*}{4}
\text{min} \quad {\mathbf{c_0^T x_0} }~~~+~~			& { c_1^T x_1 }~~+~~\cdots		&~~+~~c_N^T x_N		& \\
\text{s.t.} \quad \mathbf{A_0 x_0} ~~~~~~~~ 			&   								 	 		&										& = \mathbf{b_0} \\
\mathbf{d_0}	\leq\mathbf{C_0  \mathbf{x_0}} ~~~~~~~~ 	&  									  			&										&\leq \mathbf{f_0}  \\
{A_1  \mathbf{x_0} }~~+~~     						&  B_1 x_1 		 		        &        								&= b_1  \\
{ d_1}\leq{  C_1  \mathbf{x_0}}~~+~~& D_1 x_1 			            &        								&\leq f_1  \\
{ \vdots}~~~~~~~~~~											&~~~~~~~~~~~~~~\ddots     &   	  														&~~\vdots  \\
{ 	A_N \mathbf{x_0}} ~~+~~      					&        							            &~~+~~ B_N x_N	  		&= b_N  \\
{ 	d_N	}\leq{C_N \mathbf{x_0}}~~+~~&				                    	       	&~~+~~ D_N x_N 		&\leq f_N  \\
{ {	\mathbf{F_0 x_0}} }~~+~~     					&{  F_1 x_1}~~+~~\cdots &~~+~~F_N x_N      	&= \mathbf{b_{N+1}}  \\ 
{ \mathbf{d_{N+1}}}\leq{	\mathbf{G_0 x_0}}~~+~~   &{  G_1 x_1}~~+~~\cdots &~~+~~G_N x_N			& \leq \mathbf{f_{N+1}} \\
 && \ell_i \leq x_i &\leq u_i \quad \forall i=0,\dots,N
  \end{alignat*}
  \caption{\LP with block-diagonal structure linked by variables and constraints.}\label{fig:blockstructure}
\end{figure}
The mathematical structure of models handled by the current version of the solver are block-diagonal \LPs as specified in Fig.~\ref{fig:blockstructure}.
The $x_i \in \R^{n_i}$ are vectors of decision variables and $\ell_i$, $u_i \in (\R\cup\{\pm \infty\})^{n_i}$ are vectors of lower and upper bounds for $i = 0,1,\ldots,N$.
The extended version of \pips applies a parallel interior-point method to the problem exploiting the given structure for parallelizing expensive linear system solves.
It distributes the problem among several different processes and establishes communication between them via the Message Passing Interface (\MPI).
Distributing the \LP data among these \MPI processes as evenly as possible is an elementary feature of the solver.
Each process only knows part of the entire problem, making it possible to store and process huge \LPs that would otherwise be too large to be stored in main memory on a single desktop machine.

The \LP is distributed in the following way: For each index $i = 1,\ldots,N$ only one designated process stores the matrices $A_i$, $B_i$, $C_i$, $D_i$, $F_i$, $G_i$, the vectors $c_i$, $b_i$, $d_i$, $f_i$, and the variable bounds $\ell_i$, $u_i$.
We call such a unit of distribution a \emph{block} of the problem.
Furthermore, each process holds a copy of the block with $i=0$, containing the matrices $A_0$, $C_0$, $F_0$, $G_0$ and the corresponding vectors for bounds.
All in all, $N$ \MPI processes are used.
Blocks may be grouped to reduce $N$.
The presolving techniques presented in this paper are tailored to this special distribution and structure of the matrix.

\section{Structure-Specific Parallel Presolving}
\label{sec:presolving}

Currently, we have extended \pips by four different presolving methods.
Each incorporates one or more of the techniques described in~\cite{gurobiPres,andersen,gondzioPresolve}: singleton row elimination, bound tightening, parallel and nearly parallel row detection, and a few methods summarized under the term model cleanup.
The latter includes the detection of redundant rows as well as the elimination of negligibly small entries from the constraint matrix.

The presolving methods are executed consecutively in the order listed above.
Model cleanup is additionally called at the beginning of the presolving.
A presolving routine can apply certain reductions to the \LP: deletion of a row or column, deletion of a system entry, modification of variable bounds and the left- and right-hand side, and modification of objective function coefficients.
We distinguish between local and global reductions.
While \emph{local reductions} happen exclusively on the data of a single block, \emph{global reductions} affect more than one block and involve communication between the processes.
Since \MPI communication can be expensive, we reduced the amount of data sent and the frequency of communication to a minimum and introduced local data structures to support the synchronization whenever possible.

In the following, singleton row elimination is used as an example to outline necessary data structures and methods.
Although singleton row elimination is conceptually simple, its description still covers many difficulties arising during the implementation of preprocessing in an \HPC environment.
A singleton row refers to a row in the constraint matrix only containing one variable with nonzero coefficient.
Both for a singleton equality and a singleton inequality row, the bounds of the respective variable can be tightened.
This tightening makes the corresponding singleton row redundant and thus removable from the problem.
In the case of an equality row, the corresponding variable is fixed and removed from the system.

Checking whether a non-linking row is singleton is straightforward since a single process holds all necessary information. 
The detection of singleton linking rows requires communication between the processes.
Instead of asking all processes whether a given row is singleton, we introduced auxiliary data structures.
Let $f=(f_0,f_1,\ldots,f_N)$ denote the coefficient vector of a linking row.
Every process $i$ knows the number of nonzeros in block~$i$, i.e., $||f_i||_0$, and in block $0$, i.e., $||f_0||_0$, at all times.
At each synchronization point, every process also stores the current number of nonzeros overall blocks, $||f||_0$.
Whenever local changes in the number nonzeros of a linking row occur, the corresponding process stores these changes in a buffer, instead of directly modifying $||f_i||_0$ and $||f||_0$.
From that point on the global nonzero counters for all other processes are outdated and provide only an upper bound.
Whenever a new presolving method that makes use of these counters is entered, the accumulated changes of all processes get broadcast.
The local counters $||f_i||_0$ and $||f||_0$ are updated and stored changes are reset to zero.

\begin{figure}[t!]
	\centering
	\begin{subfigure}{0.5\textwidth}
		\begin{tikzpicture}[scale=0.72]
		\draw[blue, fill=blue!20] (-2.5,0) rectangle +(2,1);
		\draw[orange, fill=orange!20] (-2.5,1.1) rectangle +(2,1);
		\draw[green, fill=green!50] (-2.5,2.8) rectangle +(2,1);
		\draw[blue, fill=blue!20] (-2.5,4.5) rectangle +(2,1);
		\draw[green, fill=green!50] (1,0) rectangle +(1,1);
		\draw[orange, fill=orange!20] (3.5,0) rectangle +(1,1);
		\draw[green, fill=green!50] (1,2.8) rectangle +(1,1);
		\draw[orange, fill=orange!20] (3.5,1.1) rectangle +(1,1);
		
		\path (-1.5,2) -- node[auto=false]{\vdots} (-1.5,3);
		\path (-1.5,3.75) -- node[auto=false]{\vdots} (-1.5,4.75);
		\path (-0.3,0.5) -- node[auto=false]{\ldots} (0.7,0.5);
		\path (2.3,0.5) -- node[auto=false]{\ldots} (3.3,0.5);	
		\path (2.3,3) -- node[auto=false]{$\ddots$} (3.3,2);	
		\path (-0.3,4.9) -- node[auto=false]{$\ddots$} (0.7,3.5);	
		
		\draw[violet!50, line width=1mm] (-2.5,3.4)--(-0.5,3.4) (1,3.4)--(2,3.4);
		\draw[violet!50, line width=1mm] (1.5,0)--(1.5,1) (1.5,2.8)--(1.5,3.8);
		
		\draw (1.5,3.4) node[cross=2pt] {};	
		
		\node[](i) at (-2.9,3.4) {\(A_{i\cdot}\)};
		\node[](k) at (1.5,5.8) {\(A_{\cdot k}\)};
		\node[] at (-1,-.1) {};
		\end{tikzpicture}
		\caption{Singleton row leading to local changes.}
		\label{subfig:localSingRow}
	\end{subfigure}%
	~
	\begin{subfigure}{0.5\textwidth}
		\begin{tikzpicture}[scale=0.72]
		\draw[blue, fill=blue!20] (-2.5,0) rectangle +(2,1);
		\draw[orange, fill=orange!20] (-2.5,1.1) rectangle +(2,1);
		\draw[green, fill=green!50] (-2.5,2.8) rectangle +(2,1);
		\draw[blue, fill=blue!20] (-2.5,4.5) rectangle +(2,1);
		\draw[green, fill=green!50] (1,0) rectangle +(1,1);
		\draw[orange, fill=orange!20] (3.5,0) rectangle +(1,1);
		\draw[green, fill=green!50] (1,2.8) rectangle +(1,1);
		\draw[orange, fill=orange!20] (3.5,1.1) rectangle +(1,1);
		
		\path (-1.5,2) -- node[auto=false]{\vdots} (-1.5,3);
		\path (-1.5,3.75) -- node[auto=false]{\vdots} (-1.5,4.75);
		\path (-0.3,0.5) -- node[auto=false]{\ldots} (0.7,0.5);
		\path (2.3,0.5) -- node[auto=false]{\ldots} (3.3,0.5);	
		\path (2.3,3) -- node[auto=false]{$\ddots$} (3.3,2);	
		\path (-0.3,4.9) -- node[auto=false]{$\ddots$} (0.7,3.5);	
		
		\draw[violet!50, line width=1mm] (-2.5,3.4)--(-0.5,3.4) (1,3.4)--(2,3.4);
		\draw[violet!50, line width=1mm] (-1,0)--(-1,1) (-1,1.1)--(-1,2.1) (-1,2.8)--(-1,3.8) (-1,4.5)--(-1,5.5);
		
		\draw (-1,3.4) node[cross=2pt] {};	
		
		\node[](i) at (-2.9,3.4) {\(A_{i\cdot}\)};
		\node[](k) at (-1,5.8) {\(A_{\cdot k}\)};
		\node[] at (-1,-.1) {};
		\end{tikzpicture}
		\label{subfig:globalSingRow}
		\caption{Singleton row leading to global changes.}
	\end{subfigure}
	\caption{\LP data distributed on different processes and an entry \(a_{ik}\) with corresponding singleton row \(A_{i\cdot}\) and column \(A_{\cdot k}\).}\label{fig:SR_distributed}
\end{figure}
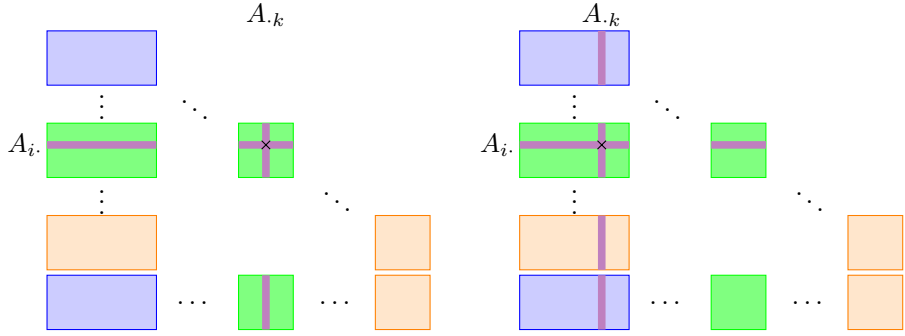

After a singleton row is detected, there are two cases to consider, both visualized in Fig.~\ref{fig:SR_distributed}.
A process might want to delete a singleton row that has its singleton entry in a non-linking part of the \LP (Fig.~\ref{fig:SR_distributed}a).
This can be done immediately since none of the other processes is affected.
By contrast, modifying the linking part of the problem is more difficult since all other processes have to be notified about the changes, e.g., when a process fixes a linking variable or when it wants to delete a singleton linking row (Fig.~\ref{fig:SR_distributed}b).
Again, communication is necessary and we implemented synchronization mechanisms for changes in variable bounds similar to the one implemented for changes in the nonzeros.

\section{Computational Results}
\label{sec:comp}

We conducted two types of experiments.
First, we looked at the general performance and impact of our presolving routines compared with the ones offered by a different \LP solver.
For the second type of experiment, we investigated the scalability of our methods.
The goal of the first experiment was to set the performance of our preprocessing into a more general context and show the efficiency of the structure-specific approach.
To this end, we compared to the sequential, source-open solver \soplex~\cite{soplex} and turned off all presolving routines that were not implemented in our preprocessing.
With our scalability experiment, we wanted to further analyze the implementation and speed-up of our presolving.
We thus ran several instances with different numbers of \MPI processes.

The instances used for the computational results come from real-world energy system models found in the literature, see~\cite{elmod1} (\texttt{elmod} instances) and~\cite{yssp} (\texttt{oms} and \texttt{yssp} instances).
All tests with our parallel presolving were conducted on the \allcaps{JUWELS} cluster at J{\"u}lich Supercomputing Centre (\allcaps{JSC}).
We used \allcaps{JUWELS}' standard compute nodes running two Intel Xeon Skylake 8168 processors each with 24 cores 2.70 \allcaps{GHz} and 96 \allcaps{GB} memory.
Since reading of the \LP and presolving it with \soplex was too time-consuming on \allcaps{JUWELS}, we had to run the tests for \soplex on a shared memory machine at Zuse institute Berlin with an Intel(R) Xeon(R) \allcaps{CPU} E7-8880 v4, 2.2\allcaps{GHz}, and 2 \allcaps{TB} of \allcaps{RAM}.

\begin{table}
	\caption{Runtimes and nonzero reductions for parallel presolving and sequential \soplex presolving. The number of nonzeros in columns ``nnzs'' are given in thousands.}\label{tab:soplexVSpips}
	\centering
	\sffamily
	\scriptsize
	\setlength{\tabcolsep}{0px}
	\begin{tabular*}{\textwidth}{@{\extracolsep{\fill}}lrrrrrrrr@{}}
		\toprule
		&\multicolumn{2}{c}{input}      &  \multicolumn{3}{c}{\pips} &\multicolumn{3}{c}{\soplex}\\
		\cmidrule(ll){2-3} \cmidrule(ll){4-6} \cmidrule(ll){7-9}
		instance & $N$ & nnzs & $t_{1}$ [s] & $t_{N}$ [s] & nnzs & $t_S$ [s] & nnzs\\
		\midrule
	oms\_1 	& 120 	& 2891k	& 1.13	& 0.02 	& 2362k	& 1.51	& 2391k	\\%& 766 & 814 & 477 & 731 & 481 & 735 & 				pips reduction: 0.817 & difference pips soplex: 0.01
	oms\_2 	& 120 	& 11432k& 5.10	& 0.19 	& 9015k 	& 11.09	& 9075k	\\%& 3174	& 3372 & 1877 & 2888 & 1892 & 2908 & 					0.788 & 0.005
	oms\_3 	& 120 	& 1696k	& 1.01	& 2.88 	& 1639k	& 0.64	& 1654k	\\%& 306 & 498 & 271 & 491 & 408 & 491 & 								0.966 & 0.009
	oms\_4 	& 120 	& 131264k& 57.25	& 3.45 	& 126242k& 206.31& 127945k\\%& 34227 & 38960 & 31495 & 38061 & 31785 & 38295 &				0.961 & 0.012
	oms\_5 	& 120 	& 216478k& 157.12& 85.41 & 158630k&$>$24h	& --	\\%& 81281 & 94965 & 37512 & 64953 & & & 							0.732 & -
	oms\_6 	& 120 	& 277923k& 187.73& 88.39 & 231796k&$>$24h	& --	\\%& 93811 & 107555 & 61317 & 77815 & & & 							0.834 & -
	elmod\_1& 438 	& 272602k& 125.62& 0.48 	& 208444k&$>$24h	& --	\\%& 98817 & 85883 & 39449 & 83773 & & & 						0.764 & -
	elmod\_2& 876 	& 716753k& 365.47& 1.05 	& 553144k&$>$24h	& --	\\%& 255916 & 226061 & 101944& 221811 & & & 					0.771 & -
	yssp\_1	& 250 	& 27927k	& 13.01 & 0.44 	& 22830k & 92.63	& 23758k	\\%& 8483 & 9478 & 5092 & 8214 & 5142	& 8348 & 			0.817 & 0.033
	yssp\_2	& 250 	& 68856k	& 33.80	& 7.28 	& 55883k	&1034.77& 59334k	\\%& 21882& 25440 & 13384 & 22167 & 13447 & 22681 & 		0.811 & 0.050
	yssp\_3	& 250 	& 32185k	& 14.10	& 0.36 	& 28874k	& 95.08 & 29802k	\\%& 9264 & 10267 & 7187 & 9301 & 7244 & 9435 & 			0.897 & 0.028
	yssp\_4	& 250 	& 85255k	& 39.71	& 7.25 	& 76504k	&1930.16& 80148k	\\%& 25005 & 28787 & 20283 & 25785 & 20376 & 26299 & 		0.897 & 0.042
	\bottomrule																													
	\end{tabular*}
\end{table}

The results of the performance experiment are shown in Table~\ref{tab:soplexVSpips}.
We compared the times spent in presolving by \soplex $t_S$, our routines running with one \MPI process $t_{1}$ and running with the maximal possible number of \MPI processes $t_{N}$.
The nnzs columns report the number of nonzeros when read in (input) and after preprocessing.
The key observations are:
\begin{itemize}
\item Except on the two smallest instances with less than 3~million nonzeros, already the sequential version of structure-specific presolving outperformed \soplex significantly.  The four largest instances with more than 200~million nonzeros could not be processed by \soplex within 24~hours.
\item The overall reduction performed by both was very similar with an average deviation of less than 2\%. The nonzero reduction overall instances was about 16\% on average.
\item Parallelization reduced presolving times on all instances except the smallest instance \texttt{oms\_3}.  On \texttt{oms\_2}, \texttt{elmod\_\{1,2\}}, and \texttt{yssp\_\{2,4\}} the speed-ups were of one order of magnitude or more.  However, on instances \texttt{oms\_\{4,5,6\}} and \texttt{yssp\_\{2,4\}} the parallel speed-up was limited, a fact that is further analyzed in the second experiment.
\end{itemize}

\begin{figure}[t]
  	\label{fig:speed-up}
	\centering
	\includegraphics[width=0.8\textwidth]{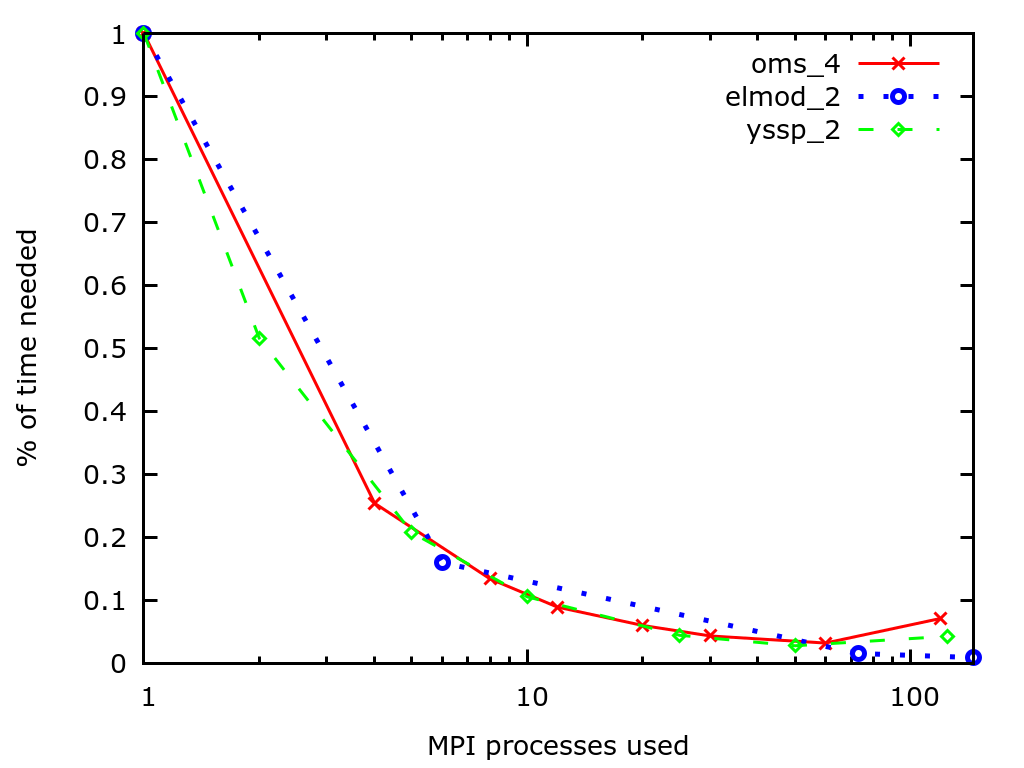}
	\caption{Total presolving time for three instances of each type, relative to time for sequential presolving with one \MPI process.}
\end{figure}

The results of our second experiment can be seen in Figure~\ref{fig:speed-up}.
We plot times for parallel presolving, normalized by the time needed by one \MPI
process.
Let $S_n = t_1/t_n$ denote the speed-up obtained with $n$ \MPI processes versus one \MPI process.
Whereas for \texttt{elmod\_2} we observe an almost linear speed-up $S_{146} \approx 114$, on \texttt{yssp\_2} and \texttt{oms\_4} the best speed-ups $S_{50}\approx 36$ and $S_{60} \approx 31$, respectively, are sublinear.
For larger numbers of \MPI processes, runtimes even start increasing again.

The limited scalability on these instances is due to a comparatively large amount of linking constraints.
As explained in Sec.~\ref{sec:presolving}, performing global reductions within linking parts of the problem increases the synchronization effort.
As a result, this phenomenon usually leads to a ``sweet spot'' for the number of
\MPI processes used, after which performance starts to deteriorate again.
This effect was also responsible for the low speed-up on \texttt{oms\_\{5,6\}} in Table~\ref{tab:soplexVSpips}.
A larger speed-up can be achieved when running with fewer processes.

To conclude, we implemented a set of highly parallel structure-preserving presolving methods that proved to be as effective as sequential variants found in an out-of-the-box \LP solver and outperformed them in terms of speed on truly large-scale problems.
Beyond the improvements of the presolving phase, we want to emphasize that the reductions helped to accelerate the subsequent interior-point code significantly.
On the instance \texttt{elmod\_1}, the interior-point time could be reduced by more than half, from about 780 to about 380~seconds.

\paragraph{Acknowledgements}
{\footnotesize This work is funded by the Federal Ministry for Economic Affairs and Energy within the \allcaps{BEAM-ME} project (ID: 03ET4023A-F) and by the Federal Ministry of Education and Research within the Research Campus \allcaps{MODAL} (\allcaps{ID}: 05M14ZAM). The authors gratefully acknowledge the Gauss Centre for Supercomputing e.V. (www.gauss-centre.eu) for funding this project by providing computing time through the John von Neumann Institute for Computing (\allcaps{NIC}) on the \allcaps{GCS} Supercomputer \allcaps{JUWELS} at J{\"u}lich Supercomputing Centre (\allcaps{JSC}).}

\ifthenelse{\zibreport = 0}{
	\bibliographystyle{splncs}
}{
	\bibliographystyle{plain}
}
\bibliography{pipspresolve}

\end{document}